\documentclass{article}

\usepackage[top=1.5in,left=1.25in,right=1.25in,bottom=1.5in]{geometry}

\usepackage{latexsym}
\usepackage[pdftex]{graphicx}
\usepackage{epsfig, ecltree, epic, eepic}
\usepackage{enumerate,amsmath,amssymb,dsfont,pifont,amsthm}
\usepackage[dvips]{color}
\usepackage{graphicx}
\usepackage[notref,notcite]{showkeys}
\usepackage[arrow, matrix, curve]{xy}
\usepackage{amsthm,amsmath}
\usepackage{stmaryrd}
\usepackage{enumerate,amssymb,dsfont,pifont}
\usepackage{xcolor}
\usepackage{tikz}
\usepackage{hyperref}
\hypersetup{
    colorlinks,
    citecolor=black,
    filecolor=black,
    linkcolor=black,
    urlcolor=black
}

\theoremstyle{plain}
\newtheorem{theorem}{Theorem}[section]

\theoremstyle{definition}
\newtheorem{definition}[theorem]{Definition}

\theoremstyle{remark}
\newtheorem{remark}[theorem]{Remark}

\newcommand{\bbp}{\mathbb{P}}

\newcommand{\bbn}{\mathbb{N}}

\makeatletter
\newcommand*\Si{\mathpalette\Si@{.5}}
\newcommand*\Si@[2]{\mathbin{\vcenter{\hbox{\scalebox{#2}{$\m@th#1\bullet$}}}}}
\makeatother

\renewcommand{\leq}{\leqslant}
\renewcommand{\geq}{\geqslant}

\begin{document}

\allowdisplaybreaks

\title{\bfseries Fixing a Minor Mistake in the Theory of Stochastic Integration and Differential Equations}

\author{%
    \textsc{Sebastian Rickelhoff}%
    \thanks{University Siegen, Department of Mathematics,
              D-57072 Siegen, Germany,
              \texttt{Sebastian.Rickelhoff@uni-siegen.de}, 
              phone: +49-231-755-3099, fax: +49-231-755-3064.} \ and
                  \textsc{Alexander Schnurr}%
    \thanks{University Siegen, Department of Mathematics,
               D-57072 Siegen, Germany,
              \texttt{Schnurr@mathematik.uni-siegen.de}, 
              phone: +49-231-755-3099, fax: +49-231-755-3064.}
    }

\date{\today}

\maketitle
%%% abstract %%%%%%%%%%%%%%%%%%%%%%%%%%%%%%%%%%%%%%%%%%%%%%%%%%%%%%%%%%%%%%%%%%%%%%%%%%
\begin{abstract}
\noindent
When considering stochastic integration and the theory of stochastic differential equations, P. Protter's textbook \cite{protter} undoubtedly is a main piece of standard literature. Not only is it well-written, but it also contains various profound results regarding these fields. Unfortunately, Theorem 12 of Chapter V, which presents an equivalence to uniform convergence on compacts in probability, is found to be incorrect. Given that numerous important results rely on this theorem, this paper aims to present a corrected version of it.
\end{abstract}


\section{Introduction}

Generations of graduate students have learned the theory of stochastic integration w.r.t. semimartingales from the well written textbook \cite{protter}. Several of them (including the authors of the present article) have used it later as reference in their mathematical research. There is probably not a single mathematical textbook without mistakes. In \cite{protter} Theorem 12 of Chapter V is not correct in its presented form (meaning in particular that the proof is also wrong). If it was only a single result - that would not be worth mentioning. However, this theorem is used to prove the Picard iteration method for SDEs (Theorem 8), various stability results (Corollary to Theorem 14, Theorem 15) and the result on the Markov nature of solutions of L\'evy driven SDEs (Theorem 32). Furthermore, the correct version of Theorem 12 has some value in its own right. 

In the present article, we first state the original result and prove by a simple example, that it is wrong. We give a hint on were the mistake in the proof is. Secondly, we state a new version of the theorem, and prove this. Since this new result is still strong enough, the proofs of the other Theorems mentioned above can be fixed easily. 

\section{The Original and Corrected Statement}
All stochastic process are considered on a stochastic basis $(\Omega, \mathcal{F}, (\mathcal{F}_t)_{t \geq 0}, \bbp)$. Let us recall the following two concepts for the reader's convenience:
\begin{definition}
\begin{itemize}
    \item[a)] A sequence $(H^n)_{n \in \bbn}$ of stochastic processes $H^n=(H_t^n)_{t\geq 0}$ is said to \textit{converge uniformly on compacts in probability (ucp)} to a process $H$ if we have 
    $$
    \lim_{n \to 0} \sup_{0 \leq s \leq t} \left| H_s^n -H_s \right| = 0, \quad \forall t >0.
    $$
    Here, the limit is a limit in probability.
    \item[b)] Let $(H_t)_{t \geq 0}$ be a càdlàg stochastic process. We define the $\underline{\underline{S}}^2$-norm of this process as follows
    $$
    \left \lVert H \right\rVert_{\underline{\underline{S}}^2} := \left\lVert \sup_{t \geq 0} \left|H_t\right| \right\rVert_{L^2}.
    $$
\end{itemize}

\end{definition}
In Chapter 2 of \cite{protter} ucp-convergence is established as \emph{the} convergence concept which is most suitable for stochastic processes in the framework of stochastic integration. The concept of convergence in $\underline{\underline{S}}^2$ is established in Chapter 5 in the context of SDEs. While ucp-convergence appears as a natural concept, the concepts in Chapter 5 are highly useful, but appear to be artificial at first glance. Theorem 12 then makes a surprising connection between the two: for $H^n, H \in \mathbb{D}$, i.e. càdlàg adapted processes $H^n$ and $H$, the author claims that ucp-convergence is equivalent to the existence of a subsequence $(n_{k_l})_{l\in \bbn}$ such that $H^{n_{k_l}}$ converges pre-locally to $H$ in $\underline{\underline{S}}^2$ for $l \to \infty$.

This statement is wrong. As a simple counterexample one may consider the process $(X_n)_{n \in \bbn}$ defined by
\[
 X_{2n}(\omega):=0 \text{ and } X_{2n+1}(\omega):=1, \quad n \in \bbn.
\]
 Obviously, the subsequence $(X_{2n})_{n\in \bbn}$ converges to the zero process in the described way while the full sequence $(X_n)_{n\in\bbn}$ does not converge ucp. 

%\section{Fixing the Mistake}

Now we state the corrected version of the theorem:

\begin{theorem}
Let $H^n, H \in \mathbb{D}$.  The following statements are equivalent: 
\begin{itemize}
\item[(i.)] $H^n \overset{\text{ucp}}\longrightarrow H$ for $n \to \infty$.
\item[(ii.)] For any subsequence $(n_k)_{k\in \bbn} \subset \bbn$ there exists a subsubsequence $(n_{k_l})_{l\in \bbn}$ such that $H^{n_{k_l}}$ converges pre-locally to $H$ in $\underline{\underline{S}}^2$ for $l \to \infty$.
\end{itemize}
\end{theorem}
\begin{proof}
We first show that `$(i.) \Rightarrow (ii.)$': Without loss of generality, we assume $H \equiv 0$. Since $H^n \overset{\text{ucp}}\longrightarrow 0$ for $n \to \infty$, $H^{n_k} \overset{\text{ucp}}\longrightarrow 0$ for $k \to \infty$ and every subsequence $(n_k)_{k\in \bbn} $. \\ 
By iteration, there exists a non-increasing sequence $(N_i)_{i \in \bbn}$ with $N_i \subset (n_k)_{k \in \bbn}$ for all $i \in \bbn$  such that 
$$
\underset{n\in N_i}{\lim_{n \to \infty}} \sup_{0\leq s \leq i} |H_s^n| = 0 \quad a.s..
$$
Due to Cantor's diagonalization procedure there exists a subsubsequence $(n_{k_l})_{l\in \bbn} \subset (n_k)_{k \in \bbn}$ such that 
$$
\lim_{l \to \infty} \sup_{0\leq s \leq i} |H_s^{n_{k_l}}| = 0 \quad a.s.,
$$
for all $i \in \bbn$.  Thus,  $(H^{n_{k_l}})_{l\in \bbn}$ converges uniformly on compacts to $0$ a.s.. Let now 
\begin{align*}
T_l &:= \inf \{ t \geq 0 : |H_t^{n_{k_l}}| \geq 1 \} \text{ and}\\
S_m &:= \inf_{j \geq m} T_j
\end{align*}
be stopping times.
Since for each $p\in \bbn$ there exists an index $N_p(\omega)$ such that 
$$
\sup_{0\leq s \leq p} \left|H_s^{n_{k_l}}(\omega)\right| < 1 \quad \forall l \geq N_p(\omega), 
$$
for almost all $\omega \in \Omega$. It follows that $T_l(\omega) \geq p$ for all $l \geq N_p(\omega)$ and, therefore, $S_{N(\omega)}(\omega) \geq p$. The sequences $(T_m)_{m\in \bbn}$ and $(S_m)_{m\in \bbn}$ both increase to $\infty$ a.s.  The sequence $(S_m \wedge m)_{m\in \bbn}$ will play the role of the localization sequence mentioned in $(ii.)$.  Therefore,  since $H^{n_{k_l}}$ converges uniformly on compacts to $0$ a.s., we observe that
$$
\sup_{s \geq 0} \left|(H_s^{n_{k_l}})^{(S_m \wedge m)-}\right|\leq  \sup_{0 \leq s \leq m}  \left|(H_s^{n_{k_l}})^{S_m}\right| \leq \sup_{0 \leq s \leq m}  \left|H_s^{n_{k_l}}\right| \underset{l \to \infty}\longrightarrow 0 \quad a.s.
$$
for all $m \in \bbn$. Moreover, for $l\geq m$ we have
\begin{align*}
\sup_{s \geq 0} \left|(H_s^{n_{k_l}})^{(S_m \wedge m)-}\right| \leq  \sup_{0 \leq s < S_m}  \left|(H_s^{n_{k_l}})^{m}\right| \leq \sup_{0 \leq s < S_m}  \left|H_s^{n_{k_l}}\right|  \leq \sup_{0 \leq s < T_l}  \left|H_s^{n_{k_l}}\right| \leq 1.
\end{align*}
Hence, Lebesgue's dominated convergence theorem provides 
\begin{align*}
\lim_{l\to \infty} \left\lVert (H^{n_{k_l}})^{(S_m \wedge m)-} \right\rVert_{\underline{\underline{S}}^2} = \lim_{l\to \infty} \left\lVert \sup_{s\geq 0} \ (H_s^{n_{k_l}})^{(S_m \wedge m)-} \right\rVert_{L^2} 
= \left\lVert\lim_{l\to \infty}  \sup_{s\geq 0} (H_s^{n_{k_l}})^{(S_m \wedge m)-} \right\rVert_{L^2} =0.
\end{align*}
This is the statement.\\
For the converse `$(ii.) \Rightarrow (i.)$', we assume that $H^n$ does not converge to $H$ in ucp, i.e. there exist $t_0>0, \ \delta_0>0, \ \varepsilon_0 >0$ and a subsequence $(n_k)_{k \in \bbn} \subset \bbn$  such that 
$$
\bbp \left( \sup_{0 \leq s \leq t_0} \left|H_s^{n_k}-H_s\right| \geq \delta_0 \right) > \varepsilon_0
$$
for all $k\in \bbn$.  By assumption, there exists a subsubsequence $(n_{k_l})_{l \in \bbn} \subset (n_k)_{k \in \bbn}$ and a sequence of stopping times $(T_m)_{m \in \bbn}$ converging to $\infty$ a.s. for $m\to \infty$ such that 
$$
\lim_{l \to \infty} \left\rVert \sup_{0 \leq s < T_m} \left|H^{n_{k_l}}_s -H_s \right| \right\lVert_{L^2} = 0
$$
for all $m \in \bbn$.
We observe that 
\begin{align*}
\varepsilon_0 &< \bbp \left( \sup_{0 \leq s \leq t_0} \left|H_s^{n_{k_l}}-H_s\right| \geq \delta_0\right)\\
&= \bbp \left( \sup_{0 \leq s \leq t_0} \left|H_s^{n_{k_l}}-H_s\right| \geq \delta_0, \ t_0 > T_m \right)+\bbp \left( \sup_{0 \leq s \leq t_0} \left|H_s^{n_{k_l}}-H_s\right| \geq \delta_0, \ t_0 \leq T_m \right)\\
&\leq \bbp (t_0 > T_m)+\bbp \left( \sup_{0 \leq s \leq T_m} \left|H_s^{n_{k_l}}-H_s \right| \geq \delta_0\right)\\
&\leq \bbp (t_0 > T_m)+ \frac{1}{\delta_0^2} \int \left(\sup_{0 \leq s \leq T_m} \left|H_s^{n_{k_l}}-H_s\right|\right)^2  d\bbp.
\end{align*}
where we used Markov's inequality in the last equation. Taking $m,l \to \infty$ provides that the last two summands converges to $0$ since $ \bbp (t_0 > T_m)$ is independent of $l$ and 
$$\int \sup_{0 \leq s \leq T_m} \left|H_s^{n_{k_l}}-H_s\right|  d\bbp$$ is independent of $m$. This is a contradiction.
\end{proof}

\begin{remark}
Except for changes made due to the correction of the original statement, the previous proof mostly follows the original one in structure and notation. However, we want to highlight the necessary changes: When proving the reverse statement \cite{protter} uses that  
$$
\int_{\{ T^k < t_0 \}} (H^{n'}-H)^*_{t_0} \ d \bbp \leq \frac{\delta \sqrt{\varepsilon}}{2},
$$
where $\bbp(\{ T^k < t_0 \}) < \frac{\delta \sqrt{\varepsilon}}{2}$. Thus, demanding $(H^{n'}-H)^*_{t_0} \leq 1$ on $\{ T^k < t_0 \}$, which, in general, is not the case.
\end{remark}
The proofs of Theorems 8, 15 and 32 as well as the Corollary to Theorem 14 of Chapter V of \cite{protter} are easily adapted.


\begin{thebibliography}{10}

\bibitem{protter}
Protter, P. \emph{Stochastic Integration and Differential Equations}. Springer, Berlin, 2nd Edition, Version 2.1,  2005.
\end{thebibliography}
\end{document}